\documentclass[11pt]{amsart}
\usepackage{amsmath,amssymb,amsfonts,url}
\usepackage{color}

\numberwithin{equation}{section}
\newtheorem{Def}{Definition}[section]
\newtheorem{thm}{Theorem}[section]
\newtheorem{lem}{Lemma}[section]

\begin{document}
\title[Monotonicity theorems on manifolds]
{Monotonicity theorems for Laplace Beltrami operator on Riemannian
manifolds}
\author{Eduardo V. Teixeira}
\address{Universidade Federal do Cear\'a\\ Departamento de
Matem\'atica\\ Av. Humberto Monte, s/n, Campus do Pici - Bloco 914
\\
Fortaleza-CE, Brazil.  CEP 60.455-760}
\email{eteixeira@pq.cnpq.br}

\author{Lei Zhang}
\address{Department of Mathematics\\
        University of Alabama at Birmingham\\
        1300 University Blvd, 452 Campbell Hall\\
        Birmingham, Alabama 35294-1170}
\email{leizhang@math.uab.edu}
\thanks{E. Teixeira acknowledges support from NSF and CNPq.
L. Zhang is supported in part by NSF Grant 0600275 (0810902)}

\date{\today}

\begin{abstract} For free boundary problems on Euclidean spaces,
the monotonicity formulas of Alt-Caffarelli-Friedman and
Caffarelli-Jerison-Kenig are cornerstones for the regularity
theory as well as the existence theory. In this article we
establish the analogs of these results for the Laplace-Beltrami
operator on Riemannian manifolds. As an application we show that
our monotonicity theorems can be employed to prove the Lipschitz
continuity for the solutions of a general class of two-phase free
boundary problems on Riemannian manifolds.
\end{abstract}

\maketitle

\section{Introduction} For two-phase free boundary problems on Euclidean spaces,
the celebrated monotonicity formula of Alt-Caffarelli-Friedman
\cite{alt} plays a fundamentally important role in the regularity
theory as well as the existence theory:

\bigskip

\noindent{\bf Theorem A } \emph{[Alt-Caffarelli-Friedman] Let
$B_1\subset \mathbb R^n$ be the unit ball, let $u_1,u_2$ be
nonnegative subharmonic functions in $C(B_1)$. Assume $u_1\cdot
u_2=0$ and $u_1(0)=u_2(0)=0$. Set
\begin{equation}\label{108e1}
\phi(r)=\frac{1}{r^4}\int_{B_r}\frac{|\nabla
u_1|^2}{|x|^{n-2}}dx\int_{B_r}\frac{|\nabla
u_2|^2}{|x|^{n-2}}dx,\quad 0<r<1.
\end{equation}
Then $\phi(r)$ is finite and and is a nondecreasing function of
$r$.}

\bigskip

There have been different extensions of this monotonicity formula
for problems with different backgrounds. For example, Caffarelli
\cite{caf1} established a monotonicity formula for variable
coefficient operators, Friedman-Liu \cite{friedman} have an
extension for eigenvalue problems. Another important extension has
been achieved by Caffarelli-Jerison-Kenig \cite{cjk} for possibly
slightly super-harmonic functions (i.e. $\Delta u_i\ge -1$, $i=1,2$)
and they derived their new form of the monotonicity theorem:

\bigskip

\noindent{\bf Theorem B }\emph{[Caffarelli-Jerison-Kenig] Suppose
the $u_1,u_2$ are non-negative, continuous functions on the unit
ball $B_1$. Suppose that $\Delta u_i\ge -1$ ($i=1,2$) in the sense
of distributions and $u_1(x)u_2(x)=0$ for all $x\in B_1$. Then there
is a dimensional constant $C$ such that
\begin{equation}\label{cjkmon}
\phi(r)\le C(n)\left(1+\int_{B_1}\frac{|\nabla
u_1|^2}{|x|^{n-2}}dx+\int_{B_1}\frac{|\nabla
u_2(x)|^2}{|x|^{n-2}}dx \right)^2,\quad 0<r\le 1.
\end{equation}
where $\phi(r)$ is defined as in (\ref{108e1}). }

\bigskip

One key estimate that the monotonicity formulas provide is the
control of $|\nabla u_1(0)|\cdot |\nabla u_2(0)|$, which is
important for the establishment of the optimal regularity results in
free boundary problems. This estimate is obtained from (\ref{108e1})
for sub-harmonic functions. However, for some real life problems
(e.g. the Prandtl-Batchelor problem
\cite{acker,batchlor1,batchlor2,caflisch,elcrat}) and some classical
problems ( e.g. see Shahgholian \cite{shahgholian}), the equations
may be inhomogeneous and we may not have $\Delta u_i\ge 0$ ($i=1,2$)
on each side of the free boundary. The ``almost monotonicity
formula" (\ref{cjkmon}) is particularly useful in these situations
and has provided a theoretical basis for the regularity theory for
many new problems (see for example \cite{cjk,cr,cs,shahgholian}).
The parabolic counterparts of (\ref{108e1}) and (\ref{cjkmon}) have
been established by Caffarelli \cite{caffarelli5}, Caffarelli-Kenig
\cite{ck} and Edquist-Petrosyan \cite{edquist} under different
contexts.

\par
It has been pointed out by Caffarelli and Salsa in \cite{salsa1}
that the tools developed for free boundary problems on Euclidean
spaces should have their counterparts for free boundary problems on
manifolds (page ix of the introduction). From theoretical and
application points of view it is natural to consider some free
boundary problems on Riemannian manifolds, rather than on Euclidean
spaces. To the best of the authors' knowledge there has not been
much progress in this direction. The purpose of this article is to
derive the analogs of the results of Alt-Caffarelli-Friedman and
Caffarelli-Jerison-Kenig for the natural operator on Riemannian
manifold: the Laplace-Beltrami operator. With the establishment of
some monotonicity formulas for the Laplace-Beltrami operator, it
becomes possible to develop the regularity theory for two-phase free
boundary problems on Riemannian manifolds.

\par
We now describe the main results of this paper. Let $(M,g)$ be a
Riemannian manifold of dimension $n\ge 2$ and let $B_1(p)$ be a
geodesic ball of radius $1$ in $M$. Let $R_m$ be the curvature
tensor. Throughout the article we use $\Lambda$ to denote the
following bound:
\begin{equation}\label{curvebound}
|R_m|+|\nabla_g R_m|\le \Lambda.
\end{equation}
For $u_1,u_2\in H^1_{loc}(B_1(p))$ we define
$$
    \phi(r)=\frac 1{r^4}\int_{B_r(p)}\frac{|\nabla_gu_1|^2}{d(x,p)^{n-2}}dV_g
    \int_{B_r(p)}\frac{|\nabla_gu_2|^2}{d(x,p)^{n-2}}dV_g
$$
where $d(x,p)$ is the geodesic distance between $x$ and $p$ under
metric $g$.
\begin{thm}\label{thm1}
Let $n\ge 2$ and $u_1,u_2\in C^0(B_1(p))$ be nonnegative functions
that satisfy
$$\Delta_gu_i\ge -1, \quad \mbox{in}\quad B_1(p), \quad i=1,2$$
in distributional sense. Suppose in addition $u_1\cdot u_2=0$, then
$u_1,u_2\in H^1_{loc}(B_1(p))$ and there exist $C(n,\Lambda)$ and
$\delta(n,\Lambda)$ such that for $0<r<\delta$,
$$\phi(r)\le C(n,\Lambda)\left(1+\int_{B_{\delta}(p)}\frac{|\nabla_gu_1|^2}{d(x,p)^{n-2}}dV_g
+\int_{B_{\delta}(p)}\frac{|\nabla_gu_2|^2}{d(x,p)^{n-2}}dV_g
\right)^2.$$
\end{thm}
\par
Our next theorem concerns sub-harmonic functions. Let $B_1(p)$ be a
geodesic ball of radius $1$ around $p$ in $(M,g)$. Let $u_1,u_2$ be
$C^0$ non-negative functions that satisfy
\begin{equation}\label{106e1}
\Delta_gu_i\ge 0,\quad \mbox{in}\quad B_1(p)
\end{equation}
in the sense of distribution. We define $\phi(r)$ as follows:
\begin{equation}\label{106e2}
\phi(r)=\frac{e^{c_0(n,\Lambda)r^2}}
{r^4}\int_{B_r(p)}\frac{|\nabla_gu_1|^2}{d(x,p)^{n-2}}dV_g \,
\int_{B_r(p)}\frac{|\nabla_gu_2|^2}{d(x,p)^{n-2}}dV_g.
\end{equation}
\begin{thm}\label{thm2}
Let $u_1,u_2\in C^0(B_1(p))$ be non-negative subharmonic functions
over $B_1(p)$ in the sense of distribution. Suppose $u_1\cdot
u_2=0$, then there exist $\delta_0(n,\Lambda)$ and $c_0(n,\Lambda)$
such that $\phi(r)$ is non-decreasing for $0<r<\delta_0$.
\end{thm}

By comparing with their Euclidean counterparts, Theorem \ref{thm2}
can be considered as an extension of the formula of
Alt-Caffarelli-Friedman, \cite{alt}, while Theorem \ref{thm1}
corresponds to the ``almost monotonicity" formula of
Caffarelli-Jerison-Kenig.

\par
As an application of Theorem \ref{thm1} we prove the Lipschitz
regularity for viscosity solutions to two-phase free boundary
problems in Riemannian manifolds. More precisely, let $(M,g)$ be a
Riemannian manifold and $\Omega$ be a bounded open set in $M$. For
a continuous function $u$ on $\Omega$, we denote
$$
    u^+(x)=\max(u(x),0);\quad
    u^-(x)=\max (-u(x),0)
$$
and
$$
    \Omega^+(u)=\{x\in \Omega;\,\, u(x)>0\};\quad
    \Omega^-(u)=\Omega\setminus \overline{\Omega^+(u)},\quad
    F(u)=\Omega\cap \partial \Omega^+(u).
$$
Given a boundary datum $h$ defined on $\partial \Omega$, bounded
functions $f_1, ~f_2 \in L^\infty(\Omega)$, and a function ${G}
\colon \mathbb{R}^2_{+} \to \mathbb{R}$, a two-phase free boundary
problem asks for a function $u \colon \Omega \to \mathbb{R}$ that
agrees with $h$ on $\partial \Omega$ and satisfies
\begin{enumerate}
    \item $\Delta_g u^{+} = f_1$ in $ \Omega^+(u)$ and $\Delta_g u^{-} = f_2$ in $
    \Omega^-(u)$.
    \item ${G} \left ( |\nabla_g u^{+}|, |\nabla_g u^{-}|
    \right ) \ge 0$ along $F(u)$.
\end{enumerate}
The equation in item $(2)$ represents the flux balance or a
transition condition from one phase to another. Notice that, due
to the prescribed flux balance $G$, $\nabla_g u$ jumps along the
free boundary, $F(u)$; therefore, Lipschitz is the optimal
regularity for an existing solution to the above free boundary
problem.
\par
Our regularity theorem for solutions to the above two-phase free
boundary problem holds in much more generality. Indeed, with the
aid of Theorem \ref{thm1}, we will show that under natural
assumptions on $G$, any weak solution to the above free boundary
problem is Lipschitz continuous (optimal regularity). The notion
we shall use for weak solutions is inspired by the viscosity
theory for free boundary problems introduced and developed by
Caffarelli in  \cite{caffarelli1, caffarelli2, caf1}.
\begin{Def} Let $\sigma:[0,\infty)\to [0,\infty)$ be such that $\sigma(\alpha)$ tends to infinity
as $\alpha$ tends to infinity. A weak solution of the free
boundary problem is a continuous function on $\Omega$ satisfying
the following:
\begin{itemize}
    \item[(a)] $-C\le \Delta_gu\le C$ as a distribution on $\Omega^+(u)$.
    \item[(b)] $-C\le \Delta_gu\le C$ as a distribution on $\Omega^-(u)$.
    \item[(c)] Let $x_0\in F(u)$. Choose $\epsilon>0$ small so that every two points in $B(x_0,\epsilon)$
     can be connected by a unique geodesic.
    Suppose there exist $\rho<\epsilon/10$ and a unit vector
    $\nu$ such that $B_{\rho}(exp_{x_0}(x_0-\rho \nu))\subset \Omega^-(u)$ and
    \begin{eqnarray*}
        u^-(x)\ge \alpha<x-x_0,\nu>^-_g+\circ(|x-x_0|)\\
        \mbox{as}\quad x\to x_0,\quad exp_{x_0}x\in B_{\rho}(exp_{x_0}(x_0-\rho \nu)),
    \end{eqnarray*}
    then
    $$
        u^+(x)\ge \sigma(\alpha)<x-x_0,\nu>^+_g+\circ(|x-x_0|).
    $$
    \item[(d)]Let $x_0\in F(u)$. Choose $\epsilon>0$ small so that every two points in
    $B(x_0,\epsilon)$ can be connected by a unique geodesic.
    Suppose there exist $\rho<\epsilon/10$ and a unit vector
    $\nu$ such that $B_{\rho}(exp_{x_0}(x_0-\rho \nu))\subset \Omega^+(u)$ and
    \begin{eqnarray*}
        u^+(x)\ge \alpha<x-x_0,\nu>^+_g+\circ(|x-x_0|),\\
        \mbox{as}\quad x\to x_0,\quad exp_{x_0}x\in B_{\rho}(exp_{x_0}(x_0+\rho \nu)),
    \end{eqnarray*}
    then
    $$
        u^-(x)\ge \sigma(\alpha)<x-x_0,\nu>^-_g+\circ(|x-x_0|).
    $$
\end{itemize}
\end{Def}

\begin{thm} \label{thm3}
If $u$ is a weak solution in $\Omega$, then $u$ is Lipschitz
continuous on any compact subset of $\Omega$.
\end{thm}

The proof of Theorem \ref{thm1} is along the same line of the proof
of Theorem 1.3 in \cite{cjk}. In their proof the argument is divided
into a ``theoretic" part and an ``arithmetic" part. The
``arithmetic" part in our case is the same so we only derive the
``theoretic" part. The difference in our situation is that we deal
with a different operator, some tools will have to be developed to
handle the difference between the Euclidean space and the manifold.
With these tools we can still fit our argument into the scheme of
\cite{cjk}.

The organization of this paper is as follows. In section two we
prove Theorem \ref{thm1}. Since the outline of our proof is similar
to that in \cite{cjk}, we shall cite the corresponding lemmas in
\cite{cjk} in the establishment of the theoretic part. In section
three we prove Theorem \ref{thm2}. The idea of using a perturbation
of $r^{-4}$ to keep the monotonicity of $\phi(r)$ was first used by
Caffarelli \cite{caf1}. In section four we prove Theorem \ref{thm3}
as an application.

\bigskip

\noindent{\bf Acknowledgement} The authors are indebted to Arshak
Petrosyan for pointing out an error in a former version of the proof
of Lemma 2.6. Petrosyan has also informed the authors about his
ongoing joint work with N. Matevosyan on monotonicity forlumas
within Euclidean spaces for operators in divergence form, which
relates to some extent to our results.

\section{Proof of Theorem \ref{thm1}}
We use the local coordinates at $p$ and treat $p$ as $0$. First we
cite the following result in \cite{hebey}:
\begin{thm}[Theorem 1.3 in \cite{hebey}] \label{Theorem 1.3 of hebey} Let $(M,g)$ be a
Riemannian manifold and $R_m$ denote its curvature tensor. Assume
$|R_m|+|\nabla_g R_m|\le \Lambda$ on the geodesic ball centered at
$0$ with radius equal to the injectivity radius at $0$. Then there
exist $K(n,\Lambda)$ and $\delta_1(n,\Lambda)$, depending only on
$n$ and $\Lambda$, such that the components $g_{ij}$ of $g$ in
geodesic normal coordinates at $0$ satisfy: For any $i,j,k=1,..,n$
and any $y\in B_{\min\left ( \delta_1,\text{inj}_{(M,g)}(0) \right
)}$ there holds
\begin{itemize}
    \item[(i)] $\frac 14\delta_{ij}\le g_{ij}(\exp_0(y))\le
    4\delta_{ij}$ (as bilinear forms).
    \item[(ii)] $ \left | g_{ij}(\exp_0(y))-\delta_{ij} \right |\le K|y|^2$  and
    $\left | \partial_kg_{ij}(\exp_0(y)) \right |\le K|y|$.
\end{itemize}
\end{thm}

Let $\delta_1(n,\Lambda)$ be the constant determined by Theorem
\ref{Theorem 1.3 of hebey}. To prove Theorem \ref{thm1}, obviously
it is enough to show that
\begin{equation}\label{929e2}
    \phi(r)\le C \left ( 1+\int_{B_{\delta_1}}\frac{|\nabla_gu_1|^2}{|x|^{n-2}}dV_g
    +\int_{B_{\delta_1}}\frac{|\nabla_gu_2|^2}{|x|^{n-2}}dV_g \right )^2,\quad
    0<r\le \delta_1.
\end{equation}
Note that in this article if we don't mention the dependence of a
constant, it is implied that this constant depends on $n$ and
$\Lambda$.
\begin{lem}\label{925lem1}
Suppose $u\in C^0(B_1)$ is non-negative and satisfies
$\Delta_gu\ge -1$ in distributional sense. Then $u\in
H^1_{loc}(B_1)$ and in distributional sense
$$2|\nabla_g u|^2\le Cu+\Delta_g(u^2) \quad B_1. $$
i.e. For $\phi\ge 0, \phi\in C^{\infty}_0(B_1)$,
\begin{equation}\label{923e1}
\int_{B_1}2|\nabla_gu|^2\phi dV_g\le \int_{B_1}Cu\phi
dV_g+\int_{B_1}u^2 \Delta_g\phi dV_g.
\end{equation}
\end{lem}

\noindent{\bf Proof of Lemma \ref{925lem1}:} To prove
(\ref{923e1}) we consider $u_m=\rho_m * u$ where
$$\rho_m(\cdot)=m^{n}\rho(m\cdot),\quad \rho\ge 0, \, \, \rho\in
C_0^{\infty}(B_1),\quad \int_{\mathbb R^n}\rho dx=1. $$  For $u_m$
we claim for $\Omega\subset\subset B_1$,
\begin{eqnarray}
u_m\to u \quad \mbox{in}\quad C^0(\Omega) \nonumber \\
\label{923e2}
 \Delta_gu_m\ge -C \quad \mbox{in}\quad \Omega \mbox{ for $m$ large}.
\end{eqnarray}

The first statement of (\ref{923e2}) is implied by the definition
of $u_m$ and the continuity of $u$. So we just derive the second
statement of (\ref{923e2}). Let $g(x)=det(g_{ij}(x))$, Theorem
\ref{Theorem 1.3 of hebey} gives the following properties
immediately: $g(x)=1+O(|x|^2)$, $\sqrt{g(x)}=1+O(|x|^2)$. Now we
have
\begin{eqnarray*}
&&\Delta_gu_m(x)=\int_{\mathbb
R^n}\frac{1}{\sqrt{g(x)}}\partial_{x_i}(\sqrt{g(x)}g^{ij}(x)\partial_{x_j}\rho_m(x-y))u(y)dy\\
&=&\int_{\mathbb R^n}\bigg
(\frac{\partial_{x_i}(\sqrt{g(x)})}{\sqrt{g(x)}}g^{ij}(x)\partial_{x_j}\rho_m(x-y)+\partial_{x_i}g^{ij}(x)
\partial_{x_j}\rho_m(x-y)\\
&&\qquad \qquad +g^{ij}(x)\partial_{x_ix_j}\rho_m(x-y)\bigg
)u(y)dy\\
&=& I_1+I_2+I_3.
\end{eqnarray*}
From the symmetry of $\rho_m$ and the continuity of $u$ we have
$I_1=\circ(1)$ and $I_2=\circ(1)$. Recall that $\circ(1)$ means a
term that tends to $0$ as $m$ tends to infinity. Thus we have
\begin{eqnarray}
\Delta_gu_m(x)&=&\int_{\mathbb
R^n}g^{ij}(x)\partial_{x_ix_j}\rho_m(x-y)u(y)dy +\circ(1)\nonumber\\
\label{924e1} &=&g^{ij}(x)\int_{\mathbb
R^n}\partial_{y_iy_j}\rho_m(x-y)u(y)dy+\circ(1).
\end{eqnarray}
On the other hand, $u$ satisfies
$$
    \int_{B_1}u\Delta_g\phi dV_g\ge -\int_{B_1}\phi dV_g, \quad
    \forall \phi \in C^{\infty}_0(B_1),
$$
which reads ($g=det(g_{ij})$)
$$
    \int_{B_1}u(y)\left (\frac{\partial_i\sqrt{g}}{\sqrt{g}}g^{ij}\partial_{y_j}\phi
    +\partial_ig^{ij}\partial_j\phi + g^{ij}\partial_{ij}\phi \right )dV_g
    \ge -\int_{B_1}\phi dV_g.
$$
By taking $\phi(y)=\rho_m(x-y)$ and letting $m$ tend to infinity
we see the first two terms are $\circ(1)$. So we have
\begin{equation}\label{924e2}
\int_{B_1}u(y)g^{ij}(y)\partial_{ij}\rho_m(x-y)dV_g\ge
-\int_{B_1}\rho_m(x-y)dV_g+\circ(1).
\end{equation}
By comparing (\ref{924e1}) and (\ref{924e2}) we obtain
\begin{equation}\label{930e1}
\Delta_gu_m(x)\ge -\sqrt{g(x)}+\circ(1).
\end{equation}
 The right hand side of
the above is $-1+O(|x|^2)+\circ(1)$. Therefore for $m$ large we
have $\Delta_gu_m(x)\ge -C$ in $\Omega$ and (\ref{923e2}) is
verified.

As a consequence of (\ref{923e2}) we have
\begin{equation}\label{924e3}
2\int_{B_1}|\nabla_g u_m|^2\phi dV_g\le \int_{B_1}(Cu_m\phi
+u_m^2\Delta_g \phi )dV_g,
\end{equation}
for $\phi\ge 0$, $\phi\in C^{\infty}_0(B_1)$. The right hand side
of (\ref{924e3}) tends to
$$\int_{B_1}(Cu\phi+u^2\Delta_g\phi
)dV_g,$$ which gives a uniform bound of the integral of
$|\nabla_gu_m|^2$ over each compact subset of $B_1$. For $u_m$ we
have
$$\int_{B_1}u_m\nabla_g\phi dV_g=-\int_{B_1}\nabla_gu_m \phi dV_g,
\quad \forall \phi\in C_0^{\infty}(B_1).$$ As a consequence, by
Riesz's representation theorem, $\nabla_gu\in L^2_{loc}(B_1)$ and
by letting $m$ tend to infinity in (\ref{924e3}) we obtain
$$\int_{B_1}|\nabla_gu|^2\phi dV_g\le \int_{B_1}(Cu\phi
+u^2\Delta_g\phi)dV_g.$$ Lemma \ref{925lem1} is established.
$\Box$

\medskip

Before stating the next lemma, we recall a standard formula, see
for instance \cite{schoenyau1}, P15.
\begin{eqnarray*}
\Delta_gr &=& \frac{n-1}r+\frac{\partial \ln
\sqrt{\det(g)}}{\partial r}\\
&=&\frac{n-1}r+O(r),\quad 0<r<\delta_1
\end{eqnarray*}
where we used $\det(g)=1+O(r^2)$. As a consequence we have
$$
    \Delta_g(r^2)=2r\Delta_gr+2|\nabla_gr|^2=2n+2r\frac{\partial \ln \sqrt{\det(g)}}{\partial r}.
$$
\begin{equation}\label{924e5}
\mbox{ For }n\ge3,\quad \Delta_g(r^{2-n})=c_n\delta_0+E_g, \quad
c_n>0.
\end{equation}
where
$$|E_g|\le Cr^{2-n}. $$
Let
\begin{equation}\label{10e1}
F_g=\left\{\begin{array}{ll}1,\qquad n=2\\
r^{2-n}+F_{1g}, \qquad n\ge 3
\end{array}
\right.
\end{equation}
where $F_{1g}$ is chosen so that
\begin{equation}\label{930e2}
F_g\ge \frac 12 r^{2-n},\quad -\Delta_gF_g\ge c_n\delta_0, \quad
\mbox{in}\quad B_{\delta_1},\quad n\ge 3.
\end{equation}
Note that for $\delta_1$ small, the estimate for $F_{1g}$ is
$$|F_{1g}|\le c_0(n,\Lambda)r^{3-n},\quad n\ge 3.$$
We shall later consider
$$
    g^t_{ij}(x)=g_{ij}(tx) \quad \text{for } |x|\le 1.
$$
The bound for $F_{1g^t}$ (for $n\ge 3$) is
\begin{equation}\label{10e2}
|F_{1g^t}(x)|\le Ct^2r^{3-n},\quad r\le 1.
\end{equation}

\begin{lem} \label{916lem1}
Let $u\in C^0(B_{\delta_1})$ be a non-negative solution of
$\Delta_gu\ge -1$ in distributional sense. Then there exist $C>0$
such that
\begin{equation}\label{916e1}
\int_{B_{\delta_1/4}}|\nabla_gu|^2F_gdV_g\le
C+C\int_{B_{\delta_1/2}\setminus B_{\delta_1/4}}u^2dV_g.
\end{equation}
\end{lem}

\noindent{\bf Proof of Lemma \ref{916lem1}:} Lemma \ref{916lem1}
corresponds to Remark 1.5 in \cite{cjk}. We claim that without
loss of generality we can assume $u$ to be smooth. Indeed if $u_m$
is the smooth approximation of $u$ considered before, we have
$u_m\to u$ a.e. and $\nabla u_m \to \nabla u$ a.e. Also $u_m$
satisfies
\begin{equation}\label{925e1}
2|\nabla_gu_m|\le 4u_m+\Delta_g(u_m^2)\quad B_{\delta_1/2}.
\end{equation}
and by (\ref{930e1})
\begin{equation}\label{925e2}
\Delta_g u_m\ge -2, \quad B_{\delta_1/4}.
\end{equation}
 Note that $u$ satisfies (\ref{925e1}) and (\ref{925e2}) in the distributional
sense. These are the inequalities we use for $u$. For $u_m$ we
shall derive
$$
    \int_{B_{\delta_1/4}}|\nabla_gu_m|^2F_gdV_g\le
    C+C\int_{B_{\delta_1/2}\setminus B_{\delta_1/4}}u_m^2dV_g.
$$
Then by letting $m\to \infty$ and applying the Dominated Convergence
Theorem we have (\ref{916e1}). Thus, hereafter in the proof of Lemma
\ref{916lem1}, we assume $u$ to be smooth.
\par
Let $\phi$ be a cut-off function such that $\phi\equiv 1$ in
$B_{\delta_1/4}$, $\phi\equiv 0$ on $B_{\delta_1/2}\setminus
B_{\frac 38\delta_1}$ and
$$|\nabla \phi|\le C,\quad |\nabla^2\phi|\le C. $$
Now we have
\begin{equation}\label{916e2}
2\int_{B_{\delta_1/4}}|\nabla_gu|^2F_g\phi dV_g\le
\int_{B_{\delta_1/4}}(4uF_g\phi+\Delta_g(u^2)\phi F_g)dV_g.
\end{equation}
To deal with the first term in the RHS of (\ref{916e2}) we shall
find a function $f$ that satisfies
$$\left\{\begin{array}{ll}
\Delta_gf\ge 2,\quad B_{\delta_1/4},\\
f(0)=0,\quad |f(x)|\le C|x|^2 \quad \mbox{in}\quad B_{\delta_1/4}.
\end{array}
\right.
$$
This function is defined as
$$f=\frac{r^2}{n-\epsilon_0}$$ where $\epsilon_0>0$ is chosen
so that
$$\Delta_gf=\frac{2n+O(r)}{n-\epsilon_0}\ge 2. $$  Consequently we have
$$\Delta_g(u+f)\ge 0 \quad \mbox{in}\quad B_{\delta_1/2}. $$
In addition to this we also have $u+f\ge 0$ in $B_{\delta_1/2}$.
With these two properties we claim that
\begin{equation}\label{917e1}
\max_{B_{\delta_1/4}}(u+f)\le C\int_{B_{\delta_1/2}\setminus
B_{\delta_1/4}}(u+f).
\end{equation}

Indeed, for $\frac 38\delta_1\le r\le \frac{\delta_1}2$, let
$f_{1,r}$ solve
$$\left\{\begin{array}{ll}
\Delta_gf_{1,r}=0,\quad B_r,\\
\\
f_{1,r}=u+f\quad \mbox{on}\quad \partial B_r.
\end{array}
\right.
$$
Then $f_{1,r}\ge u+f$ over $B_r$. Now we use the Green's
representation formula for $x\in B_{\delta_1/4}$ (see
\cite{aubin1} P112):
$$(u+f)(x)\le f_{1,r}(x)=-\int_{\partial
B_r}g^{ij}\nu_i\nabla_{jq}G(x,q)f_{1,r}(q)dS(q).$$ By the Hopf
Lemma and Theorem 4.17 of \cite{aubin1}
$$
0<-g^{ij}\nu_i\nabla_{jq}G(x,q)<C,\,\, x\in
B_{\frac{\delta_1}4},\quad \frac 38\delta_1\le |q|\le
\frac{\delta_1}2.
$$
Therefore
\begin{equation}\label{newa3}
    (u+f)(x)\le C\int_{\partial B_r}(u+f)dS.
\end{equation}
Integrating the above inequality for $\frac 38\delta_1\le r\le
\frac{\delta_1}2$ we obtain (\ref{917e1}).

With (\ref{917e1}) we go back to (\ref{916e2}) to obtain
\begin{equation}\label{newa1}
\int_{B_{\delta_1/4}}uF_g\phi\le
\max_{B_{\delta_1/4}}u\int_{B_{\delta_1/4}}F_g\phi\le
C+C\int_{B_{\delta_1/2}\setminus B_{\delta_1/4}}u^2dV_g.
\end{equation}

Next we consider $\int_{B_{\delta_1/4}}\Delta_g(u^2)\phi F_gdV_g$,
by using (\ref{930e2}) and $u(0)\ge 0$ we have
\begin{eqnarray}
\int_{B_{\delta_1/2}}\Delta_g(u^2)\phi F_gdV_g& =&
\int_{B_{\delta_1/2}}u^2\Delta_g(\phi F)dV_g \nonumber\\
&=&\int_{B_{\delta_1/2}}u^2(\Delta_g\phi F_g+2\nabla_g\phi\cdot
\nabla_gF_g+\phi\Delta_gF_g)dV_g \nonumber\\
\label{newa2}
&\le& C+C\int_{B_{\delta_1/2}\setminus
B_{\delta_1/4}}u^2.
\end{eqnarray}
Note that we used $\Delta_gF_g\le 0$ to control the last term. Lemma
\ref{916lem1} follows from (\ref{916e2}), (\ref{newa1}) and
(\ref{newa2}). $\Box$

\medskip

A consequence of Lemma \ref{925lem1} and Lemma \ref{916lem1} is
that Theorem \ref{thm1} can be proved assuming $u_1,u_2$ to be
smooth. In fact, suppose $u^i_m$ are mollified functions from
$u_i$. Then
 $\Delta u^i_m\ge -2$ over $B_{\delta_1}$. For
$u^i_m$ with $m$ large, we shall show that, for $0<r<\delta_1$,
\begin{equation}\label{929e1}
    \begin{array}{l}
        \frac 1{r^4}\displaystyle\int_{B_r}\frac{|\nabla_gu^1_m|^2}{r^{n-2}}dV_g
        \displaystyle\int_{B_r}\frac{|\nabla_gu^2_m|^2}{r^{n-2}}dV_g\\
        \le
        C \left (1+\displaystyle\int_{B_{\delta_1}}\frac{|\nabla_gu^1_m|^2}{r^{n-2}}dV_g +
        \displaystyle\int_{B_{\delta_1}}\frac{|\nabla_gu^1_m|^2}{r^{n-2}}dV_g \right
        )^2.
    \end{array}
\end{equation}

By letting $m$ tend to infinity we obtain (\ref{929e2}) from
(\ref{929e1}) by the Dominated Convergence Theorem (notice that
Lemma \ref{916lem1} makes it possible to apply the Dominated
Convergence Theorem). Thus, from now on we assume that $u_1, u_2$
are smooth positive functions which satisfy
$$
    \Delta_g u_i\ge -2 \text{ in } B_{\delta_1}.
$$
\medskip
In the remaining part of this section we shall re-scale $u_1$ and
$u_2$ as follows: For each $t<\frac{\delta_1}4$, we define
$$
    g^t_{ij}(x)=g_{ij}(tx)\quad \mbox{for }\quad |x|\le 2
$$
and
$$
    u_+(x)=u_1(tx)t^{-2},\quad u_-(x)=u_2(tx)t^{-2}, \quad |x|\le 2.
$$
In this way, we have
$$
    \Delta_{g^t}u_{\pm}(x)\ge -2,\quad x\in B_2.
$$
A key point to be noticed here is that
$$
    g^t_{ij}(x)=\delta_{ij}+O(t^2|x|^2).
$$
Since it is very cumbersome to use $g^t_{ij}$, for notational
convenience we still use $g$ to denote the metric in the remaining
part of this section, which implies that $\Delta_g$ is a
perturbation of $\Delta$ with $g_{ij}$ a perturbation of
$\delta_{ij}$, injectivity radius is greater than $4$, etc.

\begin{lem}\label{917lem1} Let $u\in W^{1,2}(B_1)$ and
$\Omega=\left \{ x\in B_1,\,\, u=0 \right \}$. Suppose
$|\Omega|\ge \mu |B_1|$ for some $\mu(n)>0$, then
$$
    \left ( \int_{B_1}|u|^pdV_g \right )^{2/p}\le C(n,\Lambda,p)\int_{\Omega}|\nabla_gu|^2dV_g,
    \quad 2\le p\le
    \frac{2n}{n-2}.
$$
\end{lem}

\noindent{\bf Proof of Lemma \ref{917lem1}:} We prove this by a
contradiction. Suppose there exists a sequence $u_k\in
W^{1,2}(B_1)$ such that $|\Omega_k|\ge \mu |B_1|$ and
$$
    \left (\int_{B_1}|u_k|^pdV_g \right )^{2/p}\ge k\int_{\Omega}|\nabla_gu_k|^2dV_g.
$$
Then let
$$
    v_k=\dfrac{u_k}{(\int_{B_1}|u_k|^pdV_g)^{1/p}}.
$$
One sees immediately that
\begin{equation}\label{917e2}
\left (\int_{B_1}|v_k|^pdV_g \right )^{1/p}=1,\quad
\int_{B_1}|\nabla_gv_k|^2dV_g\to 0.
\end{equation}
By the Sobolev-Poincar\'e inequality (see Theorem 3.7 in
\cite{hebey}):
$$
    \left (\int_{B_1}|v_k-\bar v_k|^pdv_g \right )^{1/p}\le
    C(n,\Lambda,p)\int_{B_1}|\nabla_gv_k|^2dV_g,\quad 1<p\le
    \frac{2n}{n-2}
$$
where $\bar v_k$ is the average of $v_k$ on $B_1$. Then $v_k$
converges strongly in $L^p$ norm to a constant. Since
$|\Omega_k|\ge \mu |B_1|$, this constant is $0$. However this is a
contradiction to (\ref{917e2}). Lemma \ref{917lem1} is
established. $\Box$

\begin{lem}\label{917lem2}
Let $u\in C^2(B_1)$ satisfy $\Delta_gu\ge -2$ in $B_1$. Let $F_g$
 be defined by (\ref{10e1}). Let
$\alpha=\int_{\Omega}|\nabla_gu|^2F_gdV_g<\infty$. Then there exist
$C_1,C_2$ such that if $\alpha>C_1$ and
$$\int_{\Omega\cap B_{1/4}}|\nabla_gu|^2F_gdV_g\ge \frac{\alpha}{256}$$
then $|\Omega\cap B_{1/2}\setminus B_{1/4}|>C_2$.
\end{lem}

\noindent{\bf Proof of Lemma \ref{917lem2}:} Lemma \ref{917lem2}
corresponds to Lemma 2.1 of \cite{cjk}. The proof is also similar.
We include it here for the convenience of the reader. From Lemma
\ref{916lem1} we have
$$\int_{\Omega\cap B_{\frac 14}}|\nabla_gu|^2F_gdV_g\le
C+C\int_{B_{\frac 12}\setminus B_{\frac 14}}u^2dV_g. $$ Since
$\alpha>2C$ and $\int_{B_{\frac 14}\cap
\Omega}|\nabla_gu|^2F_g>\frac{\alpha}{256}$, we have
$$\frac{\alpha}{512}\le C\int_{\Omega\cap B_{1/2}\setminus B_{1/4}}u^2dV_g. $$
If $|\Omega\cap B_{1/2}\setminus B_{1/4}|>\frac
12|B_{1/2}\setminus B_{1/4}|$, done. If not, by the Sobolev
embedding
\begin{eqnarray*}
    \frac{\alpha}{512}&\le & C \left (\int_{B_{1/2}\setminus
    B_{1/4}\cap
    \Omega}u^{\frac{2n}{n-2}}dV_g \right )^{\frac{n-2}{n}} \cdot |\Omega\cap
    B_{1/2}\setminus B_{1/4}|^{\frac 2n}\\
    &\le &C \left ( \int_{B_{1/2}\setminus B_{1/4}\cap \Omega}
    |\nabla_gu|^2dV_g \right ) \cdot |\Omega\cap B_{1/2}\setminus B_{1/4}|^{\frac
    2n}.
\end{eqnarray*}
Therefore we have $|\Omega\cap B_{1/2}\setminus B_{1/4}|>C$. Lemma
\ref{917lem2} is established. $\Box$

\begin{lem}\label{917lem3} Let $u$ be as in Lemma \ref{917lem2}.
Suppose $\alpha=\int_{B_1}|\nabla_gu|^2F_g<\infty$ and $|\Omega\cap
B_{1/2}\setminus B_{1/4}|\le (1-\lambda)|B_{1/2}\setminus B_{1/4}|$
for some $\lambda\in (0,1)$. Then there exists
$\mu(\lambda,\Lambda,n)\in (0,1)$ such that
$$\int_{\Omega\cap B_{1/4}}|\nabla_gu|^2F_gdV_g\le \mu
\int_{\Omega\cap B_{1/2}}|\nabla_gu|^2F_gdV_g.$$
\end{lem}

\noindent{\bf Proof of Lemma \ref{917lem3}:} Lemma \ref{917lem3}
corresponds to Lemma 2.3 of \cite{cjk}. Again for the convenience
of the reader we include the proof here. Since $|\Omega\cap
B_{1/2}\setminus B_{1/4}|\le (1-\lambda)|B_{1/2}\setminus
B_{1/4}|$, by Lemma \ref{917lem1} we have
$$
    \int_{\Omega\cap B_{1/2}\setminus B_{1/4}}|u|^2dV_g\le
    C_{\lambda}\int_{B_{1/2}\setminus B_{1/4}}|\nabla_gu|^2dV_g.
$$
If $\int_{B_{1/4}}|\nabla_gu|^2F_g\le \frac{\alpha}2$, there is
nothing to be proven; otherwise
$$
    \int_{B_{1/4}}|\nabla_gu|^2F_gdV_g\le C+C\int_{B_{1/2}\setminus
    B_{1/4}}u^2dV_g\le C+C\int_{B_{1/2}\setminus
    B_{1/4}}|\nabla_gu|^2dV_g
$$
Which implies
$$
    \int_{B_{1/2}\setminus B_{1/4}}|\nabla_gu|^2dV_g\ge C\alpha.
$$
Lemma \ref{917lem3} is established. $\Box$

\medskip

Let us now label important terms in our analysis:
$$
    \begin{array}{l}
        A_{\pm}(r):= \displaystyle \int_{B_r}|\nabla_g u_{\pm}|^2F_gdV_g \\
        \phi_F(r):= r^{-4}A_+(r)A_-(r), \, n\ge 3.
    \end{array}
$$

\begin{lem}\label{918lem1}
There exist constants $C_1,C_2,C_3>0$ such that if $A_{\pm}(r)\ge
C_1$ for $\frac 14\le r\le 1$, then, for a.e. $r\in (\frac 14, 1]$
\begin{equation}\label{930e8}
    \phi_F'(r)\ge -C_2\left (\frac{1}{\sqrt{A_+(r)}}+\frac{1}{\sqrt{A_-(r)}}+C_3t^2 \right)\phi_F(r).
\end{equation}
Moreover
\begin{equation}\label{930e9}
\phi_F(\frac 14)\le (1+C_2\delta )\phi_F(1)
\end{equation}
where
$\delta=\frac{1}{\sqrt{A_+(1)}}+\frac{1}{\sqrt{A_-(1)}}+C_3t^2$.
\end{lem}

\noindent{\bf Proof of Lemma \ref{918lem1}:} Lemma \ref{918lem1}
corresponds to Lemma 2.4 of \cite{cjk}. Set
$$B_{\pm}(r)=\int_{\partial B_r}|\nabla_gu_{\pm}|^2F_g\sqrt{g}dS.$$
We only consider those $r$ where $B_{\pm}(r)<\infty$. Suppose $r=1$
belongs to this set. We shall only consider $r=1$, as the estimate
for $\frac 14\le r\le 1$ is similar.
\begin{equation}\label{new2}
\phi_F'(1)=(-4)A_+A_-+B_+A_-+B_-A_+.
\end{equation}
It follows from (\ref{916e2}) that
$$
    2A_+=2\int_{B_1}|\nabla_gu_+|^2F_gdV_g\le
    \int_{B_1}\left ( 4u_+F_g+\Delta_g(u_+^2)F_g \right )dV_g.
$$
For the estimate of $A_+$ we first claim that
\begin{equation}\label{925e8}
\int_{B_1}u_+F_gdV_g\le C+C\left(\int_{\partial
B_1}u_+^2dS\right)^{\frac 12}.
\end{equation}

To see (\ref{925e8}), using the argument of Lemma \ref{916lem1} (
see (\ref{newa3})), we have
$$
    \int_{B_{3/4}}u_+F_gdV_g\le C+C \left ( \int_{\partial B_1}u_+^2dS \right )^{\frac 12}.
$$ So we only consider $x\in B_1\setminus B_{3/4}$. Let $f$ be
the function defined before so that $\Delta_g(u_++f)\ge 0$. By the
Green's representation formula for $u_++f$ we have
$$
    u_+(x)\le -\int_{\partial B_1}g^{ij}\nu_i\nabla_{jq}G(x,q)u_+(q)dS(q)+C.
$$ So
\begin{eqnarray*}
\int_{B_1\setminus B_{3/4}}u_+&\le &C+\int_{B_1\setminus B_{3/4}}
\left [-\int_{\partial
B_1}g^{ij}\nu_i\nabla_{jq}G(x,q)u_+(q)dS(q) \right ] dx\\
&\le & C-\int_{\partial B_1}u_+(q)dS(q)\int_{B_1\setminus
B_{3/4}}g^{ij}(q)\nu_i(q)\nabla_{jq}G(x,q)dx.
\end{eqnarray*}
Here we observe that
$$0< -g^{ij}(q)\nu_i(q)\nabla_{jq}G(x,q)\le C|x-q|^{1-n}.$$
This singularity makes the integral finite. This finishes the
verification of (\ref{925e8}). Therefore
\begin{eqnarray*}
    2A_+&\le &C+C \left ( \int_{\partial B_1}u_+^2dS \right )^{\frac
    12}+\int_{B_1}\Delta_g(u_+^2)F_gdV_g\\
    &\le &C+C \left ( \frac{1}{\lambda_+}\int_{\partial
    B_1}|\nabla_{\theta}u_+|^2dS \right )^{\frac 12}+T_{\lambda}
\end{eqnarray*}
where
$$
    T_{\lambda}=\int_{B_1}\Delta_g(u_+^2)F_gdV_g.
$$
Now we claim that there exists $C>0$, such that for $A_{\pm}>C$,
$\partial B_1$ meets the support of both $u_+$ and $u_-$. Indeed,
suppose without loss of generality $u_+=0$ on $\partial B_1$, since
$u_+=0$ and $\partial_{\nu}(u_+^2)=0$ on $\partial B_1$,
\begin{eqnarray*}
    2A_+&=&2\int_{B_1}|\nabla_gu_+|^2F_gdV_g \\
    &\le& \int_{B_1}(4u_++\Delta_g(u_+)^2)F_gdV_g\\
    &=&\int_{B_1}4u_+F_gdV_g+\int_{B_1}u_+^2\Delta_gF_gdV_g \\
    &\le& \int_{B_1}4u_+F_gdV_g\le C.
\end{eqnarray*}
This is a contradiction to $A_+$ being large. Let us compute
\begin{eqnarray*}
    T_{\lambda}&=& \int_{B_1}\Delta_g(u_+^2)F_gdV_g \\
    &=& \int_{B_1}\partial_i(\sqrt{g}g^{ij}\partial_j(u_+^2))F_gdx\\
    &=&\int_{\partial B_1}\sqrt{g}g^{ij}\partial_j(u_+^2)F_g\nu_idS-\int_{B_1}\sqrt{g}g^{ij}\partial_j(u_+^2)\partial_iF_gdx\\
    &=&\int_{\partial B_1}\sqrt{g}g^{ij}\partial_j(u_+^2)F_g\nu_idS
    -\int_{\partial
    B_1}\sqrt{g}(u_+^2)g^{ij}\partial_iF_g\nu_jdS+\int_{B_1}u_+^2\Delta_gF_gdV_g\\
    &\le & 2\int_{\partial
    B_1}\sqrt{g}g^{ij}u_+\partial_ju_+F_g\nu_idS-\int_{\partial
    B_1}\sqrt{g}(u_+^2)g^{ij}\partial_iF_g\nu_jdS\\
    &=& T_1+T_2
\end{eqnarray*}
For the integrant term in $T_1$ we have
\begin{eqnarray*}
    \sqrt{g}g^{ij}\partial_ju_+F_g\nu_i &=&(1+O(t^2)) \cdot (\delta_{ij}+O(t^2))
    \partial_ju_+\nu_i\\
    &=&(1+O(t^2)) \cdot (\partial_iu_+x_i+O(t^2)|\nabla u_+|)\\
    &=&\partial_ru_++O(t^2)|\nabla u_+|.
\end{eqnarray*}
Consequently
\begin{eqnarray*}
T_1&\le &2\int_{\partial B_1}u_+(\partial_ru_++O(t^2)|\nabla u|)\\
&\le &m\int_{\partial B_1}u_+^2dS+\frac 1m\int_{\partial
B_1}(\partial_ru_+)^2dS+O(t^2)\int_{\partial
B_1}u_+^2+O(t^2)\int_{\partial B_1}|\nabla u_+|^2.
\end{eqnarray*}
For $T_2$ we have
\begin{eqnarray*}
    T_2&=&-\int_{\partial B_1}\sqrt{g}(u_+^2)g^{ij}\partial_jF_g\nu_idS\\
    &=&(1+O(t^2)) \cdot (\delta_{ij}+O(t^2))\int_{\partial
    B_1}u_+^2((n-2)x_i+O(t^2))\nu_jdS\\
    &=&(n-2)\int_{\partial B_1}u_+^2dS+O(t^2)\int_{\partial
    B_1}u_+^2dS.
\end{eqnarray*}
Note that $T_2=0$ for $n=2$. Therefore
\begin{eqnarray*}
    T_{\lambda}&\le & \left [ m+n-2+O(t^2) \right ] \int_{\partial B_1}u_+^2dS+
    \left (\frac 1m+O(t^2) \right )\int_{\partial B_1}(\partial_ru_+)^2dS\\
    &+ & O(t^2)\int_{\partial B_1}|\nabla_{\theta}u|^2dS\\
    &\le &\frac{m+n-2+O(t^2)}{\lambda_+}\int_{\partial
    B_1}|\nabla_{\theta}u_+|^2dS+ \left (\frac 1m+O(t^2) \right )\int_{\partial
    B_1}(\partial_ru_+)^2dS\\
    &\le & \left ( \frac {1+O(t^2)}{\alpha_+}+O(t^2) \right )\int_{\partial B_1}|\nabla
    u_+|^2dS.
\end{eqnarray*}
where in the last step $m$ is chosen to be $m=\alpha_+$. Here we
recall that $\lambda_{\pm}=\alpha_{\pm}(\alpha_{\pm}+n-2)$.
Consequently we have
$$
    T_{\lambda}\le \left ( \frac {1+O(t^2)}{\alpha_+}+O(t^2) \right ) B_+
$$
and
$$
    2A_+\le C+C\sqrt{\frac{B_+}{\lambda_+}}+ \left ( \frac {1+O(t^2)}{\alpha_+}+O(t^2) \right )B_+.
$$
A similar estimate can also be obtained for $A_-$, so we now have
\begin{equation}\label{930e3}
2A_{\pm}\le C+C\sqrt{\frac{B_{\pm}}{\lambda_{\pm}}}+ \left (
\frac{1+O(t^2)}{\alpha_{\pm}}+O(t^2) \right )B_{\pm}.
\end{equation}
From (\ref{new2}) we see that if $B_{+}\ge 4A_{+}$ or $B_-\ge 4A_-$,
then $\phi_F'(1)\ge 0$. So we just assume $B_{\pm}\le 4A_{\pm}$.
From (\ref{930e3}) we see that if $\alpha_+\ge 3$
$$4A_+\le C+C\sqrt{\frac{A_+}{\lambda_+}}+B_+$$
for small $t$. Easy to see that Lemma \ref{918lem1} holds in this
case. The argument for $\alpha_-\ge 3$ is the same. So the only case
left is when $\alpha_{\pm}\le 3$ and $B_{\pm}\le 4A_{\pm}$. The
Friedland-Hayman inequality gives
\begin{equation}\label{926e1}
\alpha_++\alpha_-\ge 2 .
\end{equation}
 Multiply the equation
 for $A_+$ in (\ref{930e3}) by $A_-$ and multiply the equation for
 $A_-$ by $A_+$. After adding them together we have:
 \begin{eqnarray*}
    4A_+A_-&\le &2(\alpha_++\alpha_-)A_+A_-\le C+C(\sqrt{A_+}A_++\sqrt{A_-}A_+)\\
    &+& \left ( 1+O(t^2) \right )(A_-B_++B_-A_+).
\end{eqnarray*}
Using $B_{\pm}\le 4A_{\pm}$, (\ref{new2}) and (\ref{926e1}) we see
that
$$
    \phi_F'(1) \ge -C \left ( \frac 1{\sqrt{A_+}}+\frac 1{\sqrt{A_-}}+Ct^2 \right )\phi_F(1)
$$
and (\ref{930e8}) is established. From (\ref{930e8}) we divide by
$\phi_F(r)$ and integrate from $\frac 14\le r\le 1$, we have
$$\phi_F(\frac 14)\le \phi_F(1)e^{C\delta}\le
\phi_F(1)(1+C\delta)$$ because $\delta<1$. So (\ref{930e9}) holds.
Lemma \ref{918lem1} is established. $\Box$

 \medskip

 To finish the proof of Theorem \ref{thm1} we set up an iterative
 scheme as follows: Let
 $$A_k^+=\int_{|x|<4^{-k}}\frac{|\nabla_gu_1|^2}{|x|^{n-2}}dV_g,
 \quad
 A_k^-=\int_{|x|<4^{-k}}\frac{|\nabla_gu_2|^2}{|x|^{n-2}}dV_g, $$
 and $b_k^{\pm}=4^{4k}A_k^{\pm}$. The following lemma  corresponds
 Lemma 2.8 of \cite{cjk}:

 \begin{lem} \label{930lem1}
 There exist universal constants $C_1,C_2>0$ such that if $b_k^{\pm}\ge
 C_1$, then
 $$4^{4}A_{k+1}^+A_{k+1}^-\le
 A_k^+A_k^-(1+\delta_k)$$
 where
 $\delta_k=\frac{C_2}{\sqrt{b_k^+}}+\frac{C_2}{\sqrt{b_k^-}}+C_2 4^{-2k}$.
 \end{lem}

\noindent{\bf Proof of Lemma \ref{930lem1}:} Let
$$u_+(x)=4^{2k}u_1(4^{-k}x),\quad u_-(x)=4^{2k}u_2(4^{-k}x), \quad
g^k_{ij}(x)=g_{ij}(4^{-k}x). $$

As discussed before $\Delta_{g^k}u_{\pm}\ge -2$ in $B_2$. Moreover
we have
$$\int_{B_1}\frac{|\nabla_{g^k}u_+|^2}{|x|^{n-2}}dV_{g^k}=4^{4k}\int_{B_{4^{-k}}}\frac{|\nabla_gu_1|^2}{|x|^{n-2}}dV_g=b_k^+
$$ and
$$ \int_{B_1}\frac{|\nabla_{g^k}u_-|^2}{|x|^{n-2}}dV_{g^k}=4^{4k}\int_{B_{4^{-k}}}\frac{|\nabla_gu_2|^2}{|x|^{n-2}}dV_g=b_k^-.$$

By applying (\ref{930e8}) to $u_+$ and $u_-$ we have
\begin{eqnarray*}
&&4^4\int_{B_{\frac 14}}|\nabla_{g^k}u_+|^2F_{g^k}dV_{g^k}
\int_{B_{\frac 14}}|\nabla_{g^k}u_-|^2F_{g^k}dV_{g^k}\\
&\le &(1+C_n\delta_k)\int_{B_1}|\nabla_{g^k}u_+|^2F_{g^k}dV_{g^k}
\int_{B_1}|\nabla_{g^k}u_-|^2F_{g^k}dV_{g^k}
\end{eqnarray*}
 with
$$F_{g^k}=c_nr^{2-n}(1+O(4^{-2k})).$$
Therefore Lemma \ref{930lem1} is established. $\Box $

Corresponding to Lemma 2.9 in \cite{cjk} the next lemma follows
from Lemma \ref{917lem2} and Lemma \ref{917lem3} just like Lemma
2.9 of \cite{cjk} follows from Lemma 2.1 and Lemma 2.3 in that
article. So we state it without a proof.

\begin{lem}\label{930lem2}
There exist universal constants $\epsilon,C>0$ such that if
$b_k^{\pm}\ge C$ and $4^4A_{k+1}^+\ge A_k^+$, then $A_{k+1}^-\le
(1-\epsilon)A_k^-$.
\end{lem}

Theorem 1.3 in \cite{cjk} can be derived from Lemma 2.8 and Lemma
2.9 in \cite{cjk} arithmetically. The only difference here is that
there is an extra $1+O(4^{-k})$ term in Lemma \ref{930lem1}. Since
$4^{-k}$ is a geometric series, it does not affect the proof.
Theorem \ref{thm1} is established. $\Box $

\section{Monotonicity formula for sub-harmonic functions}

Since $\Delta_gu_i\ge 0$ in distributional sense, we have known from
the proof of Theorem \ref{thm1} that $u_i\in H^1_{\text{loc}}(B_1)$.
It is also easy to show that $\Delta_gu^i_m\ge -\circ(1)$, which
consequently implies
$$2|\nabla_gu_i|^2\le \Delta_g(u_i^2),\quad \mbox{in}\quad B_1$$
in weak sense. We shall assume $u_i$ be smooth in the proof of
Theorem \ref{thm2} because the argument can always be applied to
$u_m^i$ and let $m$ tend to infinity in the end.

\medskip

\subsection{Proof of Theorem \ref{thm2}}

We mainly present the proof of $n\ge 3$, as the proof of $n=2$ can
be modified easily.

From the definition of $\phi(r)$ we have
\begin{equation}\label{new3}
\phi'(r)=\bigg (-\frac
4r+2c_0r+\frac{B_r^+}{A_r^+}+\frac{B_r^-}{A_r^-}\bigg )\phi(r)
\end{equation}
where
$$A_r^+=\int_{B_r}\frac{|\nabla_gu_1|^2}{|x|^{n-2}}dV_g, \quad
A_r^-=\int_{B_r}\frac{|\nabla_gu_2|^2}{|x|^{n-2}}dV_g $$
$$B_r^+=\int_{\partial
B_r}\frac{|\nabla_gu_1|^2}{|x|^{n-2}}\sqrt{g}dS,\quad
B_r^-=\int_{\partial
B_r}\frac{|\nabla_gu_2|^2}{|x|^{n-2}}\sqrt{g}dS.
$$
We want to show that there exists $\delta_1>0$ such that
$\phi'(r)\ge 0$ for a.e. $r\in (0, \delta_1)$.

We have known that
$$-\Delta_g(|x|^{2-n})=c_n\delta_0+E$$ where
$|E|\le C|x|^{2-n}$, so we can choose $G_r=|x|^{2-n}(1+O(r^2))$ in
$B_r$ so that
$$-\Delta_gG_r\ge c_n\delta_0\quad \mbox{in}\quad B_r. $$ Now
\begin{eqnarray*}
A_r^+&\le &(1+O(r^2))\int_{B_r}|\nabla_gu_1|^2G_rdV_g\\
&\le &\frac{1+O(r^2)}2\int_{B_r}\Delta_g(u_1^2)G_rdV_g.
\end{eqnarray*}
 After
applying the integration by parts on the right hand side, we have
\begin{eqnarray*}
A_r^+&\le &\frac {1+O(r^2)}2\int_{\partial
B_r}G_r\sqrt{g}g^{ij}\partial_j(u_1^2)\nu_idS\\
&&-\frac {1+O(r^2)}2\int_{\partial
B_r}\sqrt{g}g^{ij}u_1^2\partial_iG_r\nu_jdS\\
&=&I_1+I_2.
\end{eqnarray*}
Using $g_{ij}=\delta_{ij}+O(r^2)$ and the expression of $G_r$ we
have
$$I_1=\frac 12r^{2-n}\int_{\partial
B_r}2u_1\partial_ru_1dS+O(r^{4-n})\int_{\partial B_r}u_1|\nabla
u_1|dS. $$
$$I_2=\frac{n-2}2\int_{\partial B_r}u_1^2r^{1-n}dS+O(r^{3-n}). $$
Therefore
\begin{eqnarray}
    2A_r^+&\le &(m+n-2) \left ( r^{1-n}+O(r^{3-n}) \right
    )\int_{\partial
    B_r}u_1^2+\frac 1mr^{3-n}\int_{\partial B_r}(\partial_ru_1)^2dS \nonumber\\
    \label{new1}
    &+ &O(r^{5-n})\int_{\partial B_r}|\nabla u_1|^2dS.
\end{eqnarray}
We claim that if $\alpha_+\le 3$
\begin{equation}\label{106e8}
    2\alpha_+A_r^+\le \left ( 1+O(r^2) \right )rB_r^+.
\end{equation}
To see (\ref{106e8}) holds, first we can assume $\alpha_+>0$,
otherwise it is trivial. Let $u_+(x)=u_1(rx)r^{-2}$ for $\frac
12<|x|<2$, $\bar g_{ij}(x)=g_{ij}(rx)$,
$$A_+=\int_{B_1}\frac{|\nabla_{\bar g}u_+|^2}{|x|^{n-2}}dV_{\bar
g},\quad B_+=\int_{\partial B_1}|\nabla_{\bar g}u_+|^2\sqrt{\bar
g}dS. $$ Then (\ref{106e8}) is equivalent to
\begin{equation}\label{106e9}
    2\alpha_+A_+\le \left ( 1+O(r^2) \right )B_+
\end{equation}
 where $A_+$ and $B_+$ satisfy
\begin{eqnarray*}
2A_+&\le & \left ( m+n-2+O(r^2) \right )\int_{\partial
B_1}u_+^2+\frac 1m\int_{\partial B_1}(\partial_ru_+)^2dS\\
&+&O(r^2)\int_{\partial B_1}|\nabla u_+|^2dS.
\end{eqnarray*}
Then (\ref{106e9}) can be derived by choosing $m=\alpha_+$. So
(\ref{106e8}) is established. Similarly we have
$$2\alpha_-A_r^-\le \left ( 1+O(r^2) \right )rB_r^-,\quad \mbox{if } \alpha_-\le 3. $$
If $rB_r^+>4A_r^+$ or $rB_r^->4A_r^-$ we have $\phi'(r)>0$. So we
only assume $rB_r^{\pm}\le 4A_r^{\pm}$. If $\alpha_+ \ge 3$ or
$\alpha_-\ge 3$, it is also easy to see from (\ref{new1}) that
$\phi'(r)>0$. The only case to consider is $rB_r^{\pm}\le
4A_r^{\pm}$ and $\alpha_{\pm}\le 3$. In this case, we use
(\ref{926e1}). Then, we can see that the $c_0$ term in (\ref{new3})
dominates all the $O(r^2)$ error terms. Therefore we have
$\phi'(r)>0$ for $r$ small enough. Theorem \ref{thm2} is
established. $\Box$

\section{Lipschitz continuity solutions to free boundary problems
on Riemannian manifolds}

\noindent{\bf Proof of Theorem \ref{thm3}:} This proof is a slight
modification of the proof of Theorem 4.5 in \cite{cjk}. For the
convenience of the reader we carry out all the details. Let $K$ be
a subset of $\Omega$ and let $r_0=dist(K, \partial \Omega)$. Let
$$
    K^*= \left\{ x;\,\, \text{dist}(x,K)\le \frac{r_0}2 \right \}
$$
and
$$
    M: =\max_{K^*}|u|.
$$
By the interior estimate we only need to consider points within
distance $r_0/4$ of $F(u)$. Lipschitz continuity follows from
scaled interior estimate and the bound $|u(x_1)|\le C(n,\Lambda)r$
for every $r<r_0/4$ and every $x_1\in K$ at a distance $r$ from
$F(u)$. We may assume without loss of generality that $x_1\in
\Omega^+(u)$ since the proof for $x_1\in\Omega^-(u)$ is the same.
Next, assume, for purpose of contradiction that $u(x_1)>\!>r$. It
follows from a simple scaling argument and the standard Harnack
inequality that $u(x)>\!>r$ in $B(x_1, r/2)$. Indeed, let
$$
    \bar u(y)={r^{-2}} u(ry+x_1)
$$
for $\frac 14<|y|<2$, then we see that $\bar u(0)>\!>\frac
1r>\!>C_n$. Then for $C_n$ large, it is easy to get
$$
    \bar u(y)\ge C_1(n)\bar u(0) \text{ in } B_{1/2}.
$$
Next, let $x_0\in F(u)\cap \partial B_r(x_1)$. By Hopf lemma
$$
    u^+(x)\ge \alpha<x-x_0,\nu>_g^++\circ(|x-x_0|),
$$
as $x\to x_0, \,\, x\in B_{\rho}(x_0+\rho \nu)$ for some
$\alpha>\!>1$. Therefore,
$$
    u^-(x)\ge \sigma(\alpha)<x-x_0,\nu>^+ + \circ(|x-x_0|),\quad  \sigma(\alpha)>\!>1.
$$
We use the notation $x=(x',y)\in \mathbb R^{n-1}\times \mathbb R$.
For convenience take $x_0=(0,0)$ and $\nu=(0,1)$. Let us consider
the tangent space at $x_0$ as the local coordinates. The lower bound
for $u^+$ is then
$$
    u^+(x)\ge \alpha y+O(|x|).
$$
Define
$$
    y_0(x')=\inf \left \{ y;\,\, u^+(x',y)>0,\,\, (x',y)\in B_r \right \}.
$$
Note that $\alpha>0$ and $\sigma(\alpha)>0$ imply that the graph
of $y=y_0(x')$ is tangent to $y=0$ at $(0,0)$. Since
$g_{ij}(x)=\delta_{ij}+O(|x|^2)$, we have
\begin{eqnarray*}
\int_{B_s}|\nabla_gu^+|^2dV_g\ge (1+O(s^2))\int_{B_s}|\nabla
u^+|^2dx \ge \int\limits_{B_s\cap
\{x=(x',y),y>0\}}(\alpha^2-\circ(1))dx.
\end{eqnarray*}
The last inequality above is standard: let
$l(x')=\sqrt{s^2-|x'|^2}-y_0(x')$ for $s$ small. Because $y_0(x')$
is tangent to $y=0$, $(x',y_0(x'))\in B_s$ for all
$|x'|<s-\circ(s)$. Therefore
\begin{eqnarray*}
    \int_{B_s}|\nabla u^+|^2& \ge& \int
    \int_{|x'|^2+|y|^2<s^2}|u_y^+|^2dx'dy\\
    &\ge &\int_{|x'|<s-\circ(s)}\frac 1{l(x')}\left ( \int_{y_0(x')\le y\le
    \sqrt{s^2-|x'|^2}}|u^+_y|dy \right )^2dx'\\
    &\ge &\int_{|x'|<s-\circ(s)}\frac 1{l(x')}\left (\int_{y_0(x')\le y\le
    \sqrt{s^2-|x'|^2}}u^+\right )^2dx'\\
    &\ge &\int_{|x'|<s-\circ(s)}\frac 1{l(x')} \left [\alpha
    \sqrt{s^2-|x'|^2}-\circ(s) \right ]^2dx'\\
    &=&\int\limits_{B_s\cap \{x=(x',y);y>0\}}\left [ \alpha^2-\circ(1) \right ]dx.
\end{eqnarray*}

The same bound is valid for $u^-$ with $\alpha^2$ replaced by
$\sigma(\alpha)^2$. So
\begin{eqnarray*}
\phi(R)&\ge &CR^{-4}\int_0^R\int_{B_r}|\nabla u^+|^2dx r^{1-n}dr
\,
\int_0^R\int_{B_r}|\nabla u^-|^2dx r^{1-n}dr\\
&\ge &CR^{-4}\int_0^R(\alpha^2-\circ(1))rdr \int_0^R
(\sigma(\alpha)^2-\circ(1))rdr.
\end{eqnarray*}
Thus, for sufficiently small $R$, we reach
$$
    \alpha^2\sigma^2(\alpha)\le C\phi(R).
$$
Theorem \ref{thm1} provides a uniform bound on $\phi(R)$, which
drives us to a contradiction if $\alpha$ is taken large enough.
Theorem \ref{thm3} is established. $\Box$

\end{document}